\documentclass{amsart}
\usepackage{amssymb}
\usepackage{amscd} 

\usepackage{mathrsfs}

\theoremstyle{plain}
 \newtheorem{thm}{Theorem}[section]
 \newtheorem{prop}[thm]{Proposition}
 \newtheorem{lem}[thm]{Lemma}
 
\theoremstyle{definition}
 
 \newtheorem{rem}[thm]{Remark}
 
\numberwithin{equation}{section}

\newtheorem*{thm*}{Theorem}
\newtheorem*{Hind_thm}{Hindman's Theorem}
\newtheorem*{Ellis_lem}{Ellis's Lemma}
\newtheorem{conj}[thm]{Conjecture}

\newtheorem{claim}[thm]{Claim}
\newtheorem{question}[thm]{Question}

\newcommand\one{\mathbf{1}}

\newcommand\Pscr{\mathscr{P}}

\newcommand\Acal{\mathcal{A}}
\newcommand\Bcal{\mathcal{B}}

\newcommand\Kcal{\mathcal{K}}
\newcommand\Ucal{\mathcal{U}}
\newcommand\Vcal{\mathcal{V}}

\newcommand\Abb{\mathbb{A}}
\newcommand\Fbb{\mathbb{F}}
\newcommand\Nbb{\mathbb{N}}
\newcommand\Rbb{\mathbb{R}}

\newcommand\Zbb{\mathbb{Z}}
\newcommand\Tbb{\mathbb{T}}
\newcommand\Xbb{\mathbb{X}}

\newcommand\op{\mspace{2.5mu}\widehat{\ }\mspace{2.5mu}}

\newcommand{\seq}[1]{\mathtt{#1}}

\newcommand\FIN{\operatorname{FIN}}

\newcommand\mand{\textrm{ and }}

\newcommand\lex{\mathrm{lex}}
\newcommand\num{\#}

\newcommand\Th{{}^{\mathrm{th}}}

\renewcommand{\leq}{\leqslant}
\renewcommand{\geq}{\geqslant}
\renewcommand{\setminus}{\smallsetminus}
\title[]
{Hindman's Theorem, Ellis's Lemma,\\ and Thompson's group $F$}

\author{Justin Tatch Moore}

\address{Justin Moore \\
555 Malott Hall \\
Department of Mathematics \\
Cornell University \\
Ithaca, NY 14853-4201 }

\thanks{
The research presented in this paper was partially supported by
NSF grants DMS--0757507 and DMS--1262019.
I would like to thank Stevo Todorcevic and
the math department at the University of Paris VII for their hospitality
during my visit in July 2011.
The present article underwent considerable revisions during that period.
I would also like to thank Miodrag Soki\v{c} for reading an early
draft of this paper and offering a number of suggestions.
}



\begin{document}

\maketitle

\begin{abstract}
The purpose of this article is
to formulate conjectural generalizations of Hindman's Theorem and
Ellis's Lemma for nonassociative binary systems and relate them to
the amenability problem for Thompson's group $F$.
Partial results are obtained for both conjectures.
The paper will also contain some general analysis of the conjectures.
\end{abstract}

\section{Introduction}

In \cite{amen_ramsey} a connection was established
between the amenability of discrete groups and
structural Ramsey theory.
The purpose of this article is to examine, from this perspective, the problem
of whether Richard Thompson's group $F$ is amenable.
Specifically, it will be demonstrated that the question of whether
$F$ is amenable is closely related to
the generalization of Hindman's Theorem and Ellis's Lemma to nonassociative binary systems.

Recall that a group $G$ is \emph{amenable}
if there is a finitely additive (left) translation invariant probability measure $\mu$
which measures all subsets of $G$.
Probably the most famous example of a group whose amenability is unknown
is Richard Thompson's group $F$.
The question of its amenability was considered by R. Thompson himself \cite{thompson_letter}
but was rediscovered and popularized by Geoghegan in 1979;
it first appeared in the published literature in \cite[p. 549]{comb_grp_top}.

The motivation for this question stems from the fact that $F$
does not contain a copy of $\Fbb_2$, the free group on two generators \cite{grp_pw_linear}.
It was a longstanding open problem of von Neumann to determine whether
every nonamenable group contains a copy of $\Fbb_2$
(it is easily demonstrated that any discrete group which contains $\Fbb_2$ is nonamenable).
A finitely generated counterexample was constructed by Ol$'$shanskii \cite{non_amen_nonF2} and a finitely presented
example was constructed only more recently by Ol$'$shanskii and Sapir \cite{non_amen_fp}.
Very recently, Monod constructed a new example of nonamenable group not containing
$\Fbb_2$ which is closely related to Thompson's group $F$ \cite{pw_proj}; see also \cite{fp_vN}
for a finitely presented nonamenable subgroup of Monod's group.
Thus the original motivation for considering whether $F$ is amenable
is no longer valid.
In the meantime, however, the problem of $F$'s amenability took on a life of
its own, owing to the fact that it is simple to define and serves as an important example
in group theory. 
In this paper we will see that it is related to a natural problem in Ramsey theory, giving
the problem renewed motivation.

Building on work of I. Schur and R. Rado and
confirming a conjecture of Graham and Rothschild \cite{dual_Ramsey},
Hindman proved the following result.

\begin{Hind_thm} \cite{Hindman_thm} \label{Hindman_thm}
If $c:\Nbb \to k$ is a coloring of $\Nbb$ with $k$ colors,
then there is an infinite $X \subseteq \Nbb$ such that $c$ is monochromatic on the 
sums of finite subsets of $X$.
\end{Hind_thm}

Hindman's original proof of this theorem was elementary and combinatorial but quite complex.
Galvin and Glazer later gave a simple proof using topological dynamics,
which I will now describe (see \cite[p. 102-103]{alg_betaN}, \cite[p. 30-33]{intro_Ramsey_spaces}).
The operation of addition on $\Nbb$ can be extended to its \v{C}ech-Stone compactification $\beta \Nbb$
to yield a compact left topological semigroup.
Galvin realized that the existence of an idempotent $\Ucal$
in $(\beta \Nbb,+)$ allowed for a simple recursive construction of infinite monochromatic sets as in the conclusion
of Hindman's Theorem.
Glazer then observed that the existence of such idempotents follows immediately from the following lemma
of Ellis.

\begin{Ellis_lem} \label{Ellis_lem} \cite{Ellis_lem}
If $(S,\star)$ is a compact left topological semigroup, then $S$ contains an idempotent.
\end{Ellis_lem}

We will now examine to what extent both Hindman's Theorem and Ellis's Lemma
can be generalized to a nonassociative setting.
Let $(\Tbb,\op,\one)$ denote the free binary system on one generator.
The algebra $(\Tbb,\op,\one)$ can be represented in the following manner which will be useful later
when defining our model for Thompson's group $F$.
If $a$ and $b$ are subsets of $(0,1]$,
define
\[
a \op b = \frac{1}{2} a \cup \frac{1}{2} (b + 1).
\]
Observe that, as a function, $\op$ is injective.
Consequently, the binary system generated by $\one = \{1\}$ is free;
we will take this as our model of $\Tbb$.
Just as in the case of addition on $\Nbb$, a binary operation $\star$ on a set $S$ can be
extended to $\beta S$ as follows:
$W$ is in $\Ucal \star \Vcal$ if and only if
\[
\{u \in S : \{v \in S : u \star v \in W\} \in \Vcal\} \in \Ucal
\]

Let us first observe that $(\beta \Tbb,\op)$ does not contain an idempotent
and that Hindman's theorem is false if we
replace $\Nbb$ by $\Tbb$.
To see this, define $l$ on $\Tbb$ recursively by:
\[
l(\one) = 0
\qquad
\textrm{and}
\qquad
l (a \op b) = l(a) + 1.
\]
If $\Ucal$ is an ultrafilter on $\Tbb$, then
\[
\{t \in \Tbb : l (t) \textrm{ is even}\} \in \Ucal \quad \Leftrightarrow \quad
\{t \in \Tbb : l (t) \textrm{ is odd}\} \in \Ucal \op \Ucal
\]
and in particular, $\Ucal \op \Ucal \ne \Ucal$ for any $\Ucal \in \beta \Tbb$.
Similarly, if $a,b \in \Tbb$,
then $l(a)$ and $l(a \op b)$ have different parity
and hence the naive generalization of Hindman's theorem to $(\Tbb,\op)$ fails.
Similarly, if $a,b,c \in \Tbb$,
$l((a \op b)\op c)$ and $l(a \op (b \op c))$ have a different parity.

It is informative to compare this situation to a reformulation of Hindman's Theorem.

\begin{thm} \label{Hindman_thm:Baum}
If $c : \FIN \to k$ is a coloring of $\FIN$ with $k$ colors, then
there is an infinite sequence $x_0 << x_1 << \ldots$ of elements of $\FIN$ such that
$c$ is monochromatic on all finite unions of members of this sequence.
\end{thm}
\noindent
Here $\FIN$ denotes the nonempty finite subsets of $\Nbb$ and $x << y$ abbreviates 
$\max (x) < \min (y)$.
This can be regarded as the corrected form of the following false statement:
\emph{If $c : \FIN \to k$, then there is an infinite set $X$ such that $c$ is monochromatic
on all nonempty finite subsets of $X$.}
The reason this statement is false is that
every infinite subset of $\Nbb$ contains finite nonempty subsets of both
even and odd cardinalities.
Observe that this statement is equivalent to the modification of Theorem \ref{Hindman_thm:Baum}
where we require $x_i$ to be a singleton for all $i$.
Thus we can avoid this trivial counterexample by allowing singletons to be
``glued'' together  into blocks.

For the nonassociative analog of Hindman's Theorem, I propose a different
form of ``gluing.''
Define $\Tbb_n$ to be all elements of $\Tbb$ of cardinality $n$.
These correspond to the ways to associate a sum of $n$ ones.
In particular, each $\Tbb_n$ is finite
and in fact the cardinalities of these sets are given by the Catalan numbers.
Let $\Abb_n$ denote the collection of all probability measures on $\Tbb_n$.
Notice that $\Abb_n$ can be viewed as a convex subset of the vector space generated by $\Tbb_n$
and $\Tbb_n$ can be regarded as the set of extreme points of $\Abb_n$.
In particular, if $c:\Tbb_n \to \Rbb$ is any function, then
$c$ extends linearly to a function which maps $\Abb_n$ into $\Rbb$;
such extensions will be taken without further mention.
Define $\Abb$ to be the (disjoint) union of the sets $\Abb_n$ and define
$\# : \Abb \to \Nbb$ by $\# (\nu) = n$ if $\nu \in \Abb_n$
(note that while each $\Abb_n$ is convex, $\Abb$ is not).
The operation of $\op$ on $\Tbb$ extends bilinearly to a function defined on $\Abb$:
\[
\mu \op \nu ( E ) = \sum_{a \op b \in E} \mu(\{a\}) \nu (\{b\})
\]
Observe that $\num$ is a homomorphism from $(\Abb,\op)$ to $(\Nbb,+)$.
A sequence $\mu_i$ $(i < \infty)$ of elements of $\Abb$ is
\emph{increasing} if $i < j$ implies $\#(\mu_i) < \#(\mu_j)$.

In this paper, I will prove the following partial extension
of Hindman's theorem to $(\Abb,\op)$.

\begin{thm} \label{weakNA_Hindman}
If $c:\Tbb \to [0,1]$
and $\epsilon > 0$, then there is an $r \in [0,1]$ and
an increasing sequence $\mu_i$ $(i <\infty)$ of elements of $\Abb$ such that
for all $i$,
$|c(\mu_i) - r| < \epsilon$ and all $i < j$, $|c(\mu_i \op \mu_j) - r| < \epsilon$.
\end{thm}
\noindent

In order to state the nonassociative form of Hindman's theorem
we will need the nonassociative analog of a finite sum.
If $t$ is in $\Tbb_m$, then $t$ defines a function from
$\Tbb^m \to \Tbb$ by substitution:
$t(u_0,\ldots,u_{m-1})$ is obtained by simultaneously substituting $u_i$
for the $i\Th$ occurrence of $\one$ in the term corresponding to $t$.
This operation extends to an $m$-multilinear function which maps $\Abb^m$ into $\Abb$.
We are now ready to state the conjectured generalization of Hindman's Theorem.

\begin{conj} \label{NA_Hindman}
If $c:\Tbb \to [0,1]$
and $\epsilon > 0$, then there is an $r \in [0,1]$ and an
increasing sequence $\mu_i$ $(i < \infty)$ of elements of $\Abb$ such that
whenever $t$ is in $\Tbb_m$ and
$i_0 < \ldots < i_{m-1}$ is admissible for $t$ 
\[
|c(t(\mu_{i_0},\ldots,\mu_{i_{m-1}})) - r | < \epsilon
\]
\end{conj}
Admissibility is a technical condition which will be defined later.
For now it is sufficient to mention two of its properties.
First, if $m-1 \leq i_0 < \ldots < i_{m-1}$,
then $i_0 < \ldots < i_{m-1}$ is admissible for any element of $\Tbb_m$.
Additionally any increasing sequence of $m$ integers is
admissible for some element of $\Tbb_m$.
In particular, if $c(t)$ is required to depend only on $\#(t)$,
then the above conjecture reduces to Hindman's Theorem.

It was shown in \cite{amen_ramsey} that even a weak form of Hindman's theorem
for $(\Abb,\op)$ is sufficient to prove that Thompson's group $F$ is amenable.
The relationship between $F$'s amenability and the Ramsey theory of $(\Abb,\op)$ becomes
even more apparent when one attempts to generalize Ellis's Lemma. 
The extension of a binary operation $\star$ on a set $S$ mentioned above
can be generalized so as to extend $\star$
to the space $\ell^\infty(S)^*$, which contains the set $\Pr(S)$ of
all finitely additive probability measures on $S$:
\[
\mu \star \nu (f) = \int \int f (x \star y) d \nu (y) d \mu (x)
\]
This leads to the following conjectural extension of Ellis's Lemma.
\begin{conj} \label{NA_Ellis}
If $(S,\star)$ is a binary system and $C \subseteq \Pr(S)$ is a compact convex subsystem,
then there is a $\mu$ in $C$ such that $\mu \star \mu = \mu$.
\end{conj}

The following result provides an intriguing strategy for proving $F$'s amenability.

\begin{prop}
If $\mu \in \Pr(\Tbb)$ is idempotent, then $\mu$ is $F$-invariant.
\end{prop}

The paper is organized as follows.
Section \ref{prelim_section} contains a review of the notation and background material 
which will be needed for the rest of the paper.
A proof of Theorem \ref{weakNA_Hindman} will
be given in Section \ref{weakNA_Hindman_section}.
This will serve as a warm-up for the more involved proof that
Conjecture \ref{NA_Ellis} implies Conjecture \ref{NA_Hindman} in Section \ref{implication_section}.
Section \ref{weakNA_Hindman_section} will also contain a proof
that idempotent measures in $\Pr(\Tbb)$ are $F$-invariant.
The remaining sections contain an analysis of idempotent measures and
compact convex subsystems of $\Pr(\Tbb)$.
These results fit roughly into two categories: those which offer some evidence which
makes Conjecture \ref{NA_Ellis} plausible and those which reveal what sort of difficulties
need to be addressed in proving Conjecture \ref{NA_Ellis}.

\section{Preliminaries}

\label{prelim_section}

Before beginning, let us fix some notational conventions.
In this paper, $\Nbb$ will be taken to be the positive natural
numbers and $\omega$ will denote $\Nbb \cup \{0\}$.
Elements of $\omega$ are identified with
the set of their predecessors:
$0 = \emptyset$ and $n = \{0,\ldots,n-1\}$.
If $S$ is a set, then the powerset of $S$ will be denoted by
$\Pscr (S)$.

Next we will recall the definition of a product operation on
finitely additive probability measures which is
an extension of the Fubini product of filters.
We will need some standard definitions from functional analysis;
further reading and background can be found in \cite{amenability:Paterson} \cite{functional_analysis:Rudin}.
If $X$ is a Banach space, let $X^*$ denote the collection of continuous linear functionals
on $X$.
If $S$ is a set, $\ell^\infty(S)$ denotes the space of bounded functions
from $S$ into $\Rbb$ with the supremum norm.
The space $\ell^\infty(S)^*$ will primarily be given the \emph{weak* topology}:
the weakest topology which makes the evaluation maps
$f \mapsto f(g)$ continuous for each $g$ in $\ell^\infty(S)$.
We will identify the collection $\Pr(S)$ of all finitely additive probability measures
on $S$ with the subspace
of $\ell^\infty(S)^*$ consisting of those $f$ such that
$f(g) \geq 0$ for all $g \geq 0$ and such that $f(\bar 1) = 1$, where $\bar 1$ is the function
which is constantly $1$;
if $\mu$ is in $\Pr(S)$, then $\mu$ will be identified with the bounded linear functional
$f \mapsto \int f d \mu$.
Depending on the context, we will sometimes write $f(\mu)$ for $\mu(f)$. 
The elements of $\Pr(S)$ with finite support are dense in $\Pr(S)$ in the weak* topology
and this will be
used frequently without further mention.

Suppose that $S_0$ and $S_1$ are nonempty sets.
Define $\otimes : \Pr(S_0) \times \Pr(S_1) \to \Pr(S_0 \times S_1)$
by
\[
\mu \otimes \nu (f) = \int \int f(x,y) d \nu (y) d \mu (x),
\]
where $f$ is in $\ell^\infty(S_0 \times S_1)$.
It should be noted that the order of integration is significant when measures
are required to measure all subsets of $S_0 \times S_1$.
This will be discussed further in Section \ref{monotonicity}.

\begin{prop} \label{cont_prop}
If $S_0$ and $S_1$ are nonempty sets, then
for every $\nu \in \Pr(S_1)$,
$\mu \mapsto \mu \otimes \nu$
is continuous.
Moreover if $\mu \in \Pr(S_0)$ is finitely supported,
then the map $\nu \mapsto \mu \otimes \nu$
is continuous.
\end{prop}

\begin{prop}
If $S_0$, $S_1$, and $S_2$ are nonempty sets
and $\mu_i \in \Pr(S_i)$ for $i < 3$,
then
$(\mu_0 \otimes \mu_1) \otimes \mu_2 = \mu_0 \otimes (\mu_1 \otimes \mu_2)$
up to the identification of $(S_0 \times S_1) \times S_2$ with
$S_0 \times (S_1 \times S_2)$.
\end{prop}

Now suppose that $(S,\star)$ is a binary system.
Extend $\star$ to $\Pr(S)$ as follows:
\[
\mu \star \nu (f) = \mu \otimes \nu (f \circ \star) = \mu \otimes \nu ((x,y) \mapsto f(x \star y)). 
\]
It follows from Proposition \ref{cont_prop} that if $\nu \in \Pr(S)$,
then $\mu \mapsto \mu \star \nu$ is continuous.
If $\mu$ is finitely supported, then moreover
$\nu \mapsto \mu \star \nu$ is continuous.

Thompson's group $F$ can be described as follows.
If $s,t \in \Tbb$ have equal cardinality, then
the increasing function from $s$ to $t$ extends linearly to
an automorphism of $([0,1],\leq)$.
We will write $(s \to t)$ to denote this map.
The collection of all such functions with the operation of composition is
$F$.
The group $F$ acts partially on $\Tbb$ by set-wise application with the stipulation
that $f \cdot t$
is only defined when $f$ is linear on each interval contained in the complement of $t$.
The standard generators for $F$ are given by:
\[
x_0 = \big( (\one \op \one) \op \one \to \one \op (\one \op \one)\big)
\]
\[
x_1 = \big( \one \op ((\one \op \one) \op \one) \to \one \op (\one \op (\one \op \one))\big)
\]
If we view elements of $\Tbb$ as terms,
then the partial action of $F$ on $\Tbb$ is by re-association:
\[
x_0 \cdot \big((a \op b) \op c\big) = a \op (b \op c)
\]
\[
x_1 \cdot \big(s \op ((a \op b) \op c)\big) = s \op (a \op (b \op c))
\]
The partial action of $F$ on $\Tbb$ essentially corresponds to the action of $F$ on
its positive elements with respect to the generating set $x_k$ $(k \in \omega)$, where
$x_{k+1} = x_0^k x_1 x_0^{-k}$.
It is well known that $F$ is amenable if and only if there is a $\mu$ in $\Pr(\Tbb)$
such that
\[
\mu (\{t \in \Tbb : x_0 \cdot t \textrm{ and } x_1 \cdot t \textrm{ are defined}\}) = 1,
\]
\[
\mu (x_0 \cdot E) = \mu (x_1 \cdot E) = \mu (E)
\]
whenever $E \subseteq \Tbb$ (details can be found in, e.g., \cite{fast_growth_F}).
A general introduction to $F$ and Thompson's other groups can be found in \cite{CFP}.

\section[Theorem \ref{weak_NA_Hindman}]{A partition theorem concerning sums of at most two elements of $\Abb$}

\label{weakNA_Hindman_section}

In this section we will prove Theorem \ref{weakNA_Hindman}.
The first step is to define an appropriate limit of the sequence
$\Abb_m$ $(m \in \Nbb)$.
Let $\Ucal$ be a fixed idempotent in $(\beta \Nbb,+)$.
Define $\Abb_{\Ucal}$ to be the set of all $\mu \in \Pr(\Tbb)$ such that if
$W$ is open about $\mu$, then
\[
\{m \in \Nbb : W \cap \Abb_m \ne \emptyset\} \in \Ucal.
\]
Define $\Tbb_{\Ucal} = \Abb_\Ucal \cap \beta \Tbb$.
Notice that $\Abb_{\Ucal}$ is the $\Ucal$-limit of the sequence
$\Abb_m$ $(m \in \Nbb)$ in the space
$\Kcal(\Pr(\Tbb))$ consisting of all nonempty compact subsets of $\Pr(\Tbb)$
equipped with the Vietoris topology.
In particular, $\Abb_{\Ucal}$ is compact, convex, and nonempty.

\begin{lem}
If $\mu$ and $\nu$ are in $\Abb_{\Ucal}$, so is $\mu \op \nu$.
\end{lem}

\begin{proof}
Let $W$ be open about $\mu \op \nu$ for $\mu$ and
$\nu$ in $\Abb_\Ucal$.
Let $Z = \{p \in \Nbb : W \cap \Abb_p \ne \emptyset\}$.
Since $\Ucal$ is an idempotent, we need to prove that there is a set
$X$ in $\Ucal$ such that for every $m$ in $X$,
\[
\{n \in \Nbb : m + n \in Z\} \in \Ucal
\]
Applying Proposition \ref{cont_prop}, there is an open
$U$ about $\mu$ such that if $\mu'$ is in $U$, then $\mu' \op \nu$ is in $W$.
By assumption, there is an $X$ in $\Ucal$ such that if $m$ is in $X$, then
$U \cap \Abb_m \ne \emptyset$.
For each $m$ in $X$, let $\mu_m$ be an element of $U \cap \Abb_m$.
Again by Proposition \ref{cont_prop}, there is, for each $m$ in $X$, an open $V_m$ about $\nu$ such that
if $\nu'$ is in $V_m$, then $\mu_m \op \nu'$ is in $W$.
By our assumption that $\nu$ is in $\Abb_{\Ucal}$, we have that for each $m$ in $X$ 
\[
Y_m = \{n \in \Nbb : V_m \cap \Abb_n \ne \emptyset\}
\]
is in $\Ucal$.
If $n$ is in $Y_m$, fix an element $\nu'$ of $\Abb_n \cap V_m$.
It follows that $\mu_m \op \nu'$ is in $W \cap \Abb_{m+n}$.
Thus $Y_m \subseteq \{n \in \Nbb : m + n \in Z\}$ and hence $\{n \in \Nbb : m + n \in Z\}$ is also in $\Ucal$.
\end{proof}

Theorem \ref{weakNA_Hindman} will be derived from the following proposition which is of independent
interest.

\begin{prop} \label{partial_idempotent}
If $\Bcal$ is a finite collection of subsets of $\Tbb$ and $\Ucal \in \beta \Nbb$ is idempotent,
than there is a $\mu$ in $\Abb_{\Ucal}$ such that
$\mu \op \mu \restriction \Bcal = \mu \restriction \Bcal$.
\end{prop}

\begin{proof}
By enlarging $\Bcal$ if necessary, we may assume that it is a finite Boolean subalgebra of $\Pscr(\Tbb)$.
Let $\Acal$ consist of the atoms $A$ of $\Bcal$ such that
\[
\{m \in \Nbb : A \cap \Tbb_m \ne \emptyset\}
\]
is in $\Ucal$.
Notice that if $A$ is an atom of $\Bcal$ which is not in $\Acal$, then
$\xi(A) = 0$ whenever $\xi$ is in $\Abb_{\Ucal}$.
On the other hand, if $A$ is in $\Acal$, then there is a $\xi$ in $\Tbb_{\Ucal}$ such that
$\xi(A) = 1$.
Let $X \subseteq \Tbb_{\Ucal}$ be a set of cardinality $|\Acal|$
such that for each $A$ in $\Acal$, there is a unique $\xi$ in $X$ such
that $\xi(A) = 1$.
Define a binary operation $\star : X \times X \to X$ by $\xi \star \eta = \zeta$ if $\zeta(A) =1$
where $A$ is the atom of $\Bcal$ such that $\xi \op \eta(A) = 1$.
Notice that we have $\xi \star \eta \restriction \Bcal = \xi \op \eta \restriction \Bcal$.
Extend $\star$ to a bilinear operation on the vector space generated by $X$ and let
$C$ denote the convex hull of $X$, noting that $\star$ maps $C \times C$ into $C$.
Notice that
$\mu \star \nu \restriction \Bcal = \mu \op \nu \restriction \Bcal$.
Since $X$ is finite, $\mu \mapsto \mu \star \mu$ is a continuous map defined on $C$.
Thus there is a $\mu$ in $C$ such that
$\mu \star \mu = \mu$.
It follows that $\mu \op \mu \restriction \Bcal = \mu \star \mu \restriction \Bcal = \mu \restriction \Bcal$.
\end{proof}

We are now ready to complete the proof of Theorem \ref{weakNA_Hindman}.
Let $c:\Tbb \to [0,1]$ and $\epsilon > 0$ be given.
Fix a $\delta > 0$ and a finite $\Bcal \subseteq \Pscr(\Tbb)$
such that if $\xi,\eta \in \Pr(\Tbb)$ and
$|\xi(E) - \eta(E)| < \delta$ for all $E$ in $\Bcal$,
then $|\xi(c) - \eta(c)| < \epsilon$.
If $\xi$ is in $\ell^\infty(\Tbb)^*$, let $||\xi||_{\Bcal}$ denote $\max_{E \in \Bcal} |\xi(E)|$.
Fix an element $\mu$ of $\Abb_{\Ucal}$ such that
$\mu \restriction \Bcal = \mu \op \mu \restriction \Bcal$ and set $r = c(\mu)$,
recalling the convention that
$c(\mu) = \mu(c)$.
Construct an increasing sequence $\mu_i$ $(i \in \Nbb)$ of elements of
$\Abb$ by induction so that if $i < j$, then
$||\mu(E) - \mu_i(E)||_{\Bcal} < \delta$ and
\[
||\mu \op \mu - \mu_i \op \mu ||_{\Bcal} < \delta/2 \qquad
||\mu_i \op \mu  - \mu_i \op \mu_j||_{\Bcal} < \delta/2.
\]
This is possible by Proposition \ref{cont_prop} and the definition of $\Abb_{\Ucal}$.
It follows that if $i \in \Nbb$, then
$||\mu_i - \mu||_{\Bcal} < \delta$ and hence $|c(\mu_i) - r| < \epsilon$.
Simiarly if $i < j$, then
\[
||\mu_i \op \mu_j - \mu||_{\Bcal} \leq ||\mu - \mu_i \op \mu ||_{\Bcal} + ||\mu_i \op \mu - \mu_i \op \mu_j||_{\Bcal} < \delta
\]
and consequently
$|c(\mu_i \op \mu_j) - r| < \epsilon$.
This finishes the proof of Theorem \ref{weakNA_Hindman}.

Extending this result even to right (or left) associated sums of at most three elements of $\Abb$
seems to require new ideas.
In particular, it is not clear how to find, for a given $c \in \ell^\infty(\Tbb)$, a $\mu$ in $\Abb_{\Ucal}$
such that
$c(\mu)  = c(\mu \op \mu ) = c(\mu \op (\mu \op \mu) )$.

I will conclude this section with the following proposition which
relates the nonassociative form of Ellis's Lemma
to the amenability problem for $F$.

\begin{prop} \label{idempotent_inv}
If $\mu \in \Pr(\Tbb)$ is an idempotent measure, then $\mu$ is $F$-invariant.
\end{prop}

\begin{proof}
Suppose that $\mu \in \Pr(\Tbb)$ satisfies $\mu \op \mu = \mu$;
we need to show that $\mu$ is $F$-invariant.
First observe that
\[
\mu( \{\one\} ) = \mu \op \mu (\{\one\}) = \mu \otimes \mu (\emptyset) = 0.
\]
Also
\[
\mu \op \mu (\{t \in \Tbb : \exists a (t = a \op \one)\})
=\mu (\Tbb) \cdot \mu (\{\one\}) = 1 \cdot 0 = 0.
\]
Since $\mu$ is an idempotent, following identities hold:
\[
\mu = \mu \op (\mu \op \mu) = (\mu \op \mu) \op \mu
\]
\[
\mu = \mu \op (\mu \op (\mu \op \mu)) =
\mu \op ((\mu \op \mu) \op \mu)
\]
Now suppose that $E \subseteq \Tbb$.
\[
\mu (E) = \mu (
\{t \in E : \exists a \exists b \exists c (t = (a \op b) \op c)\}
)
\]
\[
= (\mu \op \mu) \op \mu (
\{t \in E : \exists a \exists b \exists c (t = (a \op b) \op c)\}
)
\]
\[
= (\mu \otimes \mu) \otimes \mu (
\{(a,b,c) \in \Tbb^3 : (a \op b) \op c \in E\})
\]
\[
=  \mu \otimes (\mu \otimes \mu) (
\{(a,b,c) \in \Tbb^3 : a \op (b \op c) \in x_0 \cdot E\})
\]
\[
= \mu \op (\mu \op \mu) (\{t \in x_0 \cdot E : \exists a \exists b \exists c (t  = a \op (b \op c))\})
=
\mu ( x_0 \cdot E).
\]
A similar computation 
shows that $\mu(x_1 \cdot E) = \mu(E)$.
\end{proof}

\section{Conjecture \ref{NA_Ellis} implies Conjecture \ref{NA_Hindman}}

\label{implication_section}

In this section I will prove that
the nonassociative form of Ellis's Lemma implies the nonassociative
form of Hindman's Theorem.
This will be an elaboration on the proof of Theorem \ref{weakNA_Hindman} from
Proposition \ref{partial_idempotent}.

Before proceeding, it is necessary to define the notion of \emph{admissibility} from
the statement of Conjecture \ref{NA_Hindman}.
In order to motivate the definition of admissibility, consider the natural adaptation of the
proof in the previous section:
at stage $n$, we have constructed measures $\mu_i$ $(i < n)$ and wish to pick a measure $\mu_n$
in $\Abb$ such that if $t$ is any element of $\Tbb$ and $i_k$ $(k < m)$ is an increasing sequence
of length at most $\#(t)$, then
$c(t(\mu_{i_0},\ldots,\mu_{i_{m-1}},\mu,\mu,\ldots,\mu))$
and $c(t(\mu_{i_0},\ldots,\mu_{i_{m-1}},\mu_n,\mu, \ldots,\mu)$ differ by less than
$\epsilon 2^{-n-1}$.
The problem is that there are infinitely many $t$'s to consider.
The fix to this problem is to consider only those $t$ of cardinality at most $n$ at stage $n$.
This readily allows us to prove the form of Conjecture \ref{NA_Hindman}
in which one obtains the conclusion for those sequences $i_k$ $(k < m)$ satisfying $m-1 \leq i_0$.

It is possible to do better, however, by noticing that if $\mu$ is idempotent, then expressions such
as $t(\mu_{i_0},\ldots,\mu_{i_{m-1}},\mu,\mu,\ldots,\mu)$ may allow for some algebraic simplification.
This is made precise as follows.
For each $k < \#(t)$, define $t \mapsto t_k$ by
\[
t_k =
\begin{cases}
\one & \textrm{ if } t = \one \textrm{ and } k = 0 \\
a_k \op \one & \textrm{ if } t = a \op b \textrm{ and } k < \#(a) \\
a \op b_{k-\#(a)}     & \textrm{ if } t = a \op b  \textrm{ and } \#(a) \leq k < \#(t) \\
\end{cases}
\]
Notice that if we regard $t$ as a rooted orderd binary tree,
then $t_k$ is the result of iteratively removing all carets in $t$ which involve only leaves
of index greater than $k$ (where leaves are indexed in increasing order from left to
right).
Set $l_k(t) = \#(t_k) - 2$. 
If $t$ is in $\Tbb_m$,
then an increasing sequence $i_k$ $(k < m)$ is \emph{admissible} for $t$ if
for all $k < m$, $l_k(t) \leq i_k$.

Notice that $l_k(t) < \#(t)$ for each $k$ and in particular a sequence
$i_k$ $(k < m)$ is admissible for any element of $\Tbb_m$ provided that
$m \leq i_0$.
Also, the value of $l_k$ at any right associated power of $\one$ is at most $k$
and hence any increasing sequence is admissible for some element of $\Tbb_m$.

It will be helpful to adopt the following notation:
if $t$ is in $\Tbb_m$, $k \leq m$, and $\nu_i$ $(i < k)$ and $\mu$ are in $\Pr(\Tbb)$,
let $t(\nu_0,\ldots,\nu_{k-1};\mu)$ denote
$t(\nu_0,\ldots,\nu_{k-1},\mu,\ldots,\mu)$
(i.e. the sequence $\nu_i$ $(i < k)$ is extended to a sequence of length $m$ by adding
on a sequence of $m-k$ many $\mu$'s and then substituting into $t$).
The key property of the definition of $t_k$ is that 
whenever $t$ in $\Tbb$, $\nu_i$ $(i < k)$ are in $\Abb$, and $\mu$ is an idempotent in $\Pr(\Tbb)$,
then
$t(\nu_0,\ldots,\nu_{k-1};\mu) = t_k(\nu_0,\ldots,\nu_{k-1};\mu)$;
this is easily established by induction on $\#(t)$.

We are now ready to prove the main result of the section.

\begin{thm} \label{Ellis_to_Hindman}
Conjecture \ref{NA_Ellis} implies Conjecture \ref{NA_Hindman}.
\end{thm}

\begin{proof}
In Section \ref{weakNA_Hindman_section}, we proved that if $\Ucal$ is an idempotent in $(\beta \Nbb,+)$, then 
$\Abb_{\Ucal}$ is a nonempty compact convex subsystem of $(\Pr(\Tbb),\op)$.
Thus if Conjecture \ref{NA_Ellis} is true, then $\Abb_\Ucal$ contains a $\mu$
such that $\mu \op \mu = \mu$.
Fix a $c:\Tbb \to [0,1]$ and set $r = c(\mu)$.
Construct an increasing sequence $\mu_i$ $(i \in \omega)$
in $\Abb$ by recursion such that, if $\mu_i$ $(i < n)$ have been constructed, then
for all $k < m \leq n+2$ and $i_0 < \ldots < i_{k-1} < n$ and $t$ in $\Tbb_m$,
\[
|c(t(\mu_{i_0},\ldots,\mu_{i_{k-1}};\mu)) - c(t(\mu_{i_0},\ldots,\mu_{i_{k-1}},\mu_n;\mu))| < \epsilon 2^{-n-1}.
\]
This is possible by applying the definition of $\Abb_{\Ucal}$ and
the following claim.

\begin{claim} \label{cont_claim}
If $t$ is in $\Tbb_m$ and
$\nu_i$ $(i < m)$ are such that $\nu_i$ is in $\Pr(\Tbb)$ and
has finite support if $i < k$, then the function $F$ defined by
\[
F(\zeta) = t(\nu_0,\ldots,\nu_{k-2},\zeta,\nu_k,\ldots,\nu_{m-1})
\]
is continuous.
\end{claim}

\begin{proof}
The proof is by induction on $m$.
If $m = 1$, then there is nothing to show since then $F$
is just the identity.
If $m > 1$ and $t$ is in $\Tbb_m$,
then there are $a$ and $b$ such that $t = a \op b$.
If $\#(a) = l \leq k$, then
\[
t(\nu_0,\ldots,\nu_{k-2},\zeta,\nu_k,\ldots) = 
a(\nu_0,\ldots,\nu_{l-1}) \op b(\nu_l,\ldots,\nu_{k-1},\zeta,\nu_k,\ldots)
\]
Letting $\nu = a(\nu_0,\ldots,\nu_{l-1})$, we have that
\[
F(\zeta) = \nu \op b(\nu_l,\ldots,\nu_{k-1},\zeta,\nu_k,\ldots,\nu_{m-1}).
\]
Continuity of $F$ now follows from Proposition \ref{cont_prop}
and the induction hypothesis applied to $b$.
A similar argument handles the case $\#(a) > k$.
\end{proof}

Now we will verify that $\mu_k$ $(k \in \omega)$ satisfies the conclusion of
Conjecture \ref{NA_Hindman}.
To this end, let $t$ be an element of $\Tbb_m$ and
let $i_0 < \ldots < i_{m-1}$ be admissible for $t$.
By construction we have 
\[
|c(t(\mu_{i_0},\ldots,\mu_{i_{m-1}})) - c(t(\mu,\ldots,\mu))| \leq
\]
\[
\sum_{k < m} |c(t(\mu_{i_0},\ldots,\mu_{i_{k}};\mu))-
c(t(\mu_{i_0},\ldots,\mu_{i_{k-1}};\mu))|
\]
\[
= \sum_{k < m} |c(t_k(\mu_{i_0},\ldots,\mu_{i_{k}};\mu))-
c(t_k(\mu_{i_0},\ldots,\mu_{i_{k-1}};\mu))|. 
\]
By admissibility, $\#(t_k) = l_k(t) + 2 \leq i_k +2$ and
thus $\mu_{i_k}$ was chosen such that
\[
|c(t_k(\mu_{i_0},\ldots,\mu_{i_k};\mu))-
c(t_k(\mu_{i_0},\ldots,\mu_{i_{k-1}};\mu))| < \epsilon 2^{-i_k -1}.
\]
Recalling that $r = c(t(\mu,\ldots,\mu))$ and putting this all together we have that
\[
|c(t(\mu_{i_0},\ldots,\mu_{i_{m-1}})) - r| < \epsilon \sum_{k} 2^{-k-1} = \epsilon.
\]
\end{proof}

We finish this section by recalling a result of \cite{amen_ramsey}
which asserts that a weak form of the finitary version of
Conjecture \ref{NA_Hindman}
is in fact equivalent to the amenability of $F$.
If $m \in \Nbb$, an \emph{embedding of $\Tbb_m$ in $\Tbb$} is
function of the form
\[
t \mapsto t(u_0,\ldots,u_{m-1})
\]
for some sequence $u_i$ $(i < m)$ of elements of $\Tbb$.
Observe that such an embedding maps into $\Tbb_n$ for some $n$.
An embedding of $\Tbb_m$ into $\Abb_n$ is a convex combination of embeddings
of $\Tbb_m$ into $\Tbb_n$.
The range of such an embedding is a \emph{copy of $\Tbb_m$ in $\Abb_n$}.

\begin{thm} \cite{amen_ramsey}
The following are equivalent:
\begin{enumerate}

\item \label{F_amen}
Thompson's group $F$ is amenable.

\item \label{Tm_ramsey}
For every $m$ there is an $n$ such that
if $c:\Tbb_n \to [0,1]$, then there is a copy $\Xbb$
of $\Tbb_m$ in $\Abb_n$ such that
$|c(\nu) - c(\nu')| \leq 1/2$
whenever $\nu,\nu' \in \Xbb$.

\item For every $m$ there is an $n$ such that
if $c:\Tbb_n \to \{0,1\}$, then there is a copy $\Xbb$
of $\Tbb_m$ in $\Abb_n$ such that $c$ is constant on $\Xbb$.

\end{enumerate}
\end{thm}

\begin{rem}
A more restrictive notion of copy is the following:
a \emph{strong copy} of $\Tbb_m$ in $\Abb$ is the range of a function of the form
\[
t \mapsto t(\mu_0,\ldots,\mu_{m-1})
\]
where $\mu_i$ $(i < m)$ are elements of $\Abb$.
The finite form of Conjecture \ref{NA_Hindman} would
assert that for every $m$ 
there is an $n$ such that if $f:\Tbb_n \to [0,1]$,
then there is a strong copy of $\Tbb_m$ in $\Abb_n$ on which
$f$ is within $1/2$ of being constant.
It is unclear if this assertion is equivalent to the amenability of $F$.
\end{rem}

\section{Monotonicity properties of $F$-invariant measures}

\label{monotonicity}

In \cite{fast_growth_F}, a lower bound was established for the F\o lner function
for $F$, assuming that it is amenable.
This was achieved by establishing and analyzing
the following qualitative property of $F$-invariant measures in $\Pr(\Tbb)$, which is of
independent interest.

\begin{prop} \cite{fast_growth_F}
If $\mu$ is an $F$-invariant measure in $\Pr(\Tbb)$, then for $\mu$-a.e.
$t \in \Tbb$, one of the following pairs of inequalities holds:
\[
\#(t/\seq{001}) < \#(t/\seq{01}) < \#(t/\seq{10})
\]
\[
\#(t/\seq{10}) < \#(t/\seq{01}) < \#(t/\seq{001}).
\]
\end{prop}

Here $t/\sigma$ is defined recursively as follows:
$t/\sigma = t$ if $\sigma$ is the empty string and
\[
(a \op b)/\seq{0 \sigma} = a/\sigma \qquad (a \op b)/\seq{1 \sigma} = b/\sigma.
\]
That is, $t/\sigma$ is the sub-term of $t$ located at address $\sigma$.
The proof of the above proposition generalizes to yield the following result.
Recall that a \emph{quasi-order} is a reflexive, transitive relation; a quasi-order $\preceq$ is \emph{linear}
if for every $x$ and $y$, either $x \preceq y$ or $y \preceq x$.

\begin{prop} \label{Finv_monotone}
Suppose that $\preceq$ is a linear quasi-order on $\Tbb$ such that
if $s,t \in \Tbb$, then $s,t \prec s \op t$.
If $\mu \in \Pr(\Tbb)$ is an $F$-invariant measure, then either
\begin{itemize}

\item
for every incompatible pair $\sigma <_\lex \varsigma$ of nonconstant finite binary sequences,
$\mu$-a.e. $t$ satisfy
$t/\sigma \prec t/\varsigma$ or

\item
for every incompatible pair $\sigma <_\lex \varsigma$ of nonconstant finite binary sequences,
$\mu$-a.e. $t$ satisfy $t/\varsigma \prec t/\sigma$.

\end{itemize}
\end{prop}

\begin{rem}
It is worth noting an example of a linear quasi-order on $\Tbb$ which is quite different
than the order defined by $\#(s) \leq \#(t)$.
Define an equivalence relation $\equiv_{LD}$ on $\Tbb$ by relating two elements if they can be proved equal using
the \emph{left self distributive law} $a \op (b \op c) = (a \op b) \op (a \op c)$.
Define $a <_{LD} b$ if there are $a'$ and $b'$ which are
$\equiv_{LD}$-equivalent to $a$ and $b$, respectively, such
that $a'$ is a subterm of $b'$ (i.e. $a' = b'/\sigma$ for some $\sigma$ of positive length).
By work of Laver and Dehornoy \cite{alg_elem_embedd} \cite{braid_grp_LD},
if $a$ and $b$ are in $\Tbb$, then exactly one of the following is true:
$a <_{LD} b$, $a \equiv_{LD} b$, or $b <_{LD} a$.
Set $a \leq_{LD} b$ if $a \equiv_{LD} b$ or $a <_{LD} b$.
It follows that $\leq_{LD}$ a linear quasi-order which satisfies
the hypothesis of Proposition \ref{Finv_monotone}.
Notice that except for the left associated elements of $\Tbb$, which do not allow any application of
the left self distributive law, the equivalence classes of $\equiv_{LD}$ are infinite.
\end{rem}

The purpose of this section is to prove the following proposition, which
shows that a natural strategy for proving Conjecture \ref{NA_Ellis}
does at least yield measures with this qualitative property.

\begin{prop}
Suppose $C \subseteq \Pr(\Tbb)$ is minimal with respect to the properties
of being compact, convex, $\op$-closed, and nonempty and let
$\preceq$ be a linear quasi-order on $\Tbb$ such that for all $s,t \in \Tbb$,
$s,t \prec s \op t$.
Either
\begin{itemize}

\item
for every $\mu \in C$ and incompatible pair $\sigma <_\lex \varsigma$ of finite binary sequences,
$\mu$-a.e. $t$ satisfy
$t/\sigma \prec t/\varsigma$ or

\item
for every $\mu \in C$ and incompatible pair $\sigma <_\lex \varsigma$ of finite binary sequences,
$\mu$-a.e. $t$ satisfy $t/\varsigma \prec t/\sigma$.

\end{itemize}
\end{prop}

\begin{proof}
Let $C$ be given as in the statement of the proposition.
First observe that the convex hull of $C \op C$ is dense in $C$.
This is because the closure of the convex hull of $C \op C$ is clearly
nonempty, compact, convex, and $\op$-closed
(in fact if $C \op C \subseteq X \subseteq C$, then $X$ is $\op$-closed).
We will need the following claim.

\begin{claim} \label{0-1}
If $I \subseteq \Tbb$ is an interval in $(\Tbb,\preceq)$, then
either $\mu(I) = 1$ for every $\mu$ in $C$ or $\mu(I) = 0$ for every $\mu$ in $C$.
\end{claim}

\begin{proof}
First observe that it suffices to prove the claim if $I$ is an initial segment of $\Tbb$ in the $\preceq$-order.
Suppose that there exists a $\mu$ in $C$ such that $\mu(I) = p < 1$.
Let $C' = \{\nu \in C : \nu(I) \leq p^2\}$.
Observe that $C'$ is compact and convex.
Next notice that if $s$ and $t$ are in $\Tbb$ and $s \op t$ is in $I$,
then both $s$ and $t$ are in $I$.
Thus
\[
\mu \op \mu (I) = \mu \op \mu (\{s \op t \in \Tbb : s \op t \in I\})
\]
\[
\leq \mu \op \mu (\{s \op t \in \Tbb : s,t \in I\}) = \mu \op \mu (I \op I) = \mu(I) \cdot \mu(I) = p^2
\]
and hence $\mu \op \mu$ is in $C'$.
Furthermore, if $\xi$ and $\eta$ are in $C'$, then a similar computation shows that
$\xi \op \eta(I) \leq p^4 \leq p^2$ and thus that $C'$ is $\op$-closed.
We have therefore showed that $C' \subseteq C$ is compact, nonempty, and $\op$-closed.
It follows that $C' = C$ and thus that $\mu$ is in $C'$.
But this means $p \leq p^2$.
Since $p < 1$ holds by assumption, we have that $p=0$.
It follows that every element of $C = C'$ assigns measure $0$ to $I$.
\end{proof}

Let $I$ denote the set of all $s$ in $\Tbb$ such that
for every $\mu$ in $C$,
\[
\mu(\{t \in \Tbb : s \prec t\}) = 1
\]
and observe that $I$ is an initial part of $\Tbb$.
Define $J = \Tbb \setminus I$.

The remainder of the proof breaks into two cases.
First consider the case that $\mu(I) = 1$ for every $\mu$ in $C$.
We will need the following claim.
\begin{claim}
For each $\mu$ in $C$,
each $s$ in $I$, and
each finite binary sequence $\sigma$
\[
\mu(\{t \in \Tbb : t/\sigma \in I \mand s \prec t/\sigma\}) = 1.
\]
\end{claim}

\begin{proof}
Fix $s \in I$;
the proof is by induction on the length of $\sigma$.
We have already established the base case.
In order to handle the inductive step,
let $\sigma$ be a given finite binary sequence of positive length.
Notice that the set $A$ of $\mu$ in $C$ such that
\[
\mu(\{t \in \Tbb : t/\sigma \in I \mand s \prec t/\sigma\}) = 1
\]
is compact and convex.
Thus it suffices to show that $C \op C \subseteq A$.
Observe that if $\sigma = \seq{0} \bar \sigma$, then
$(u \op v)/\sigma = u/\bar \sigma$ and if $\sigma = \seq{1} \bar \sigma$,
then $(u \op v)/\sigma = v/\bar \sigma$.
If $\xi \op \eta$ is in $C \op C$ and $\sigma = \seq{0} \bar \sigma$,
then the claim follows from
\[
\xi \op \eta(\{t \in \Tbb : t/\sigma \in I \mand s \prec t/\sigma\}) =
\xi (\{u \in \Tbb : u/\bar \sigma \in I \mand s \prec u/\bar \sigma\})
\]
which equals $1$ by the induction hypothesis.
The case $\sigma = \seq{1} \bar \sigma$ is analogous.
\end{proof}
We will now prove the first alternative of the proposition by induction
on the maximum length of the sequences.
As in the claim, the set of $\mu$ which satisfy the first alternative of the 
proposition is compact and convex and hence it is sufficient to show that it
contains $C \op C$.
Let $\sigma <_\lex \varsigma$ be given finite binary sequences.
If $\sigma$ and $\varsigma$ both begin with the same digit,
then observe that if $\xi \op \eta$ is in $C \op C$, then the induction hypothesis
implies that for $\xi \op \eta$-a.e. $t$, $t/\sigma \prec t/\varsigma$.
If $\sigma = \seq{0} \bar \sigma$ and $\varsigma = \seq{1} \bar \varsigma$,
then fix a $u$ such that $u/\bar \sigma$ is in $I$.
If $\xi \op \eta$ is in $C \op C$, then
for $\eta$-a.e. $v$, $u/\bar \sigma \prec v/\bar \varsigma$.
Notice that if $t = u \op v$, then $t/\sigma = u/\bar \sigma$ and
$t/\varsigma = v/\bar \varsigma$.
Since $\xi$-a.e. $u$ satisfies that $u/\bar \sigma$ is in $I$,
it follows that for $\xi \op \eta$-a.e. $t$,
$t/\sigma \prec t/\varsigma$.
This finishes the proof of the case $\mu(I) = 1$ for
all $\mu$ in $C$.

Now suppose that $\mu(J) = 1$ for all $\mu$ in $C$.
This case is almost identical to the proof of the previous case, but with
the observation that, by definition of $I$ and Claim \ref{0-1}, we have that
$\mu ( \{t \in \Tbb : s \preceq t\}) = 0$ whenever $s$ is in $J$ and $\mu$ is in $C$.
\end{proof}

\section{The set of idempotents in $\Abb_{\Ucal}$}

The purpose of this section is to show that if $\Abb_{\Ucal}$ contains any idempotents,
it contains at least $2^{\aleph_0}$ of them.
This is then used to show that if there any idempotents in $\Abb_{\Ucal}$,
then the set of idempotents in $\Abb_{\Ucal}$ is far from being closed.

\begin{prop}
Suppose that $\{\mu_i : i \in \omega \}$ consists of idempotents and
that there are sets $E_i \subseteq \Tbb$ such that
$\mu_i(E_j) = 1$ if $i=j$ and $0$ otherwise.
Any limit point of $\{\mu_i : i \in \omega\}$ is not an idempotent.
\end{prop}

\begin{proof}
Suppose that $\mu_i$ $(i \in \omega)$ and $E_i$ $(i \in \omega)$ are 
as in the statement of the proposition and let $\mu$ be a limit point
of $\{\mu_i : i \in \omega\}$.
Let $t_i$ $(i \in \omega)$ list $\Tbb$ and define
\[
E = \bigcup_{j < i} t_i \op E_j.
\]
Observe that for each $i$, $\mu (\bigcup_{j < i} E_j) = 0$ and therefore
$\mu \op \mu (E) = 0$.
On the other hand, for each $k \in \omega$, $\mu_k( \bigcup_{j < i} E_j) = 1$ whenever
$k < i$.
Since idempotent measures assign measure $0$ to finite sets, it follows that
\[
\mu_k(E) = \mu_k \op \mu_k(E) = \int \int \chi_E(s \op t)\ d \mu_k(t)\ d \mu_k(s) = 1.
\]
Since $\mu$ is a limit point of $\{\mu_k : k \in \omega\}$, it follows
that $\mu(E) = 1 \ne \mu \op \mu(E) = 0$.
\end{proof}

\begin{prop}
If $\Abb_{\Ucal}$ contains an idempotent, then there exist idempotents
$\{\mu_r : r \in 2^{\omega}\} \subseteq \Abb_{\Ucal}$ and sets $E_r \subseteq \Tbb$
for $r \in 2^{\omega}$ such that $\mu_r(E_s) = 1$ if $r = s$ and $0$ otherwise.
In particular, the set of idempotents in $\Abb_{\Ucal}$ is either empty or else is
not a closed set.
\end{prop}

\begin{proof}
Begin by defining $u_\sigma$ by induction on the length of $\sigma$,
whenever $\sigma$ a nonempty finite binary sequence.
If $\sigma$ has length $1$, define $u_\sigma = \one$.
If $u_\sigma$ has been defined, set
\[
u_{\seq{0} \sigma} = u_\sigma \op \one \qquad u_{\seq{1} \sigma} = \one \op u_\sigma
\]
and observe that $\#(u_\sigma)$ coincides with the length of $\sigma$.
For $r$ in $2^{\omega}$ and $n \in \omega$, define $E_{r,n}$ to be the subsystem
of $(\Tbb,\op)$ generated by $\{u_{r \restriction k} : n < k\}$ and
let $E_r$ denote $E_{r,0}$.
Observe that if $r \ne s$, then for sufficiently large $n$, $E_{r,n}$ is disjoint from $E_s$.

For $r \in 2^{\omega}$, define $h_r: \Tbb \to E_r$ by
\[
h_r(t) =
\begin{cases}
h_r(a) \op h_r(b) & \textrm{ if } t = a \op b \textrm{ and } a << b \\
u_{r \restriction \#(t)} & \textrm{ if } t \textrm{ is not as above.}
\end{cases} 
\]
Here $a << b$ means that for some $p$, $\#(a) < 2^p$ and $2^p$ divides $\#(b)$.
Extend $h_r$ to a continuous linear operator on $\ell^\infty(\Tbb)$,
also denoted by $h_r$ as follows:
\[
h_r(\mu)(f) = \mu (f \circ h_r) = \int f(h_r(t))\ d \mu (t) 
.
\]
Notice that $h_r$ maps each $\Abb_m$ into $\Abb_m$ and hence maps
$\Abb_{\Ucal}$ into $\Abb_{\Ucal}$.

The proposition follows from the next two claims.
First observe that if $p$ is in $\Nbb$ and $\mu$ is in $\Abb_{\Ucal}$,
then
\[
\mu (\{t \in \Tbb : 2^p \textrm{ divides } \#(t)\}) = 1 
\]
This follows from the fact that $\Zbb/2^p\Zbb$ has a unique idempotent and
that the canonical homomorphism from $(\Nbb,+)$ into
$\Zbb/2^p\Zbb$ extends to a homomorphism of $(\beta\Nbb,+)$ into
$\Zbb/2^p \Zbb$ which must send $\Ucal$ to $0$.
 
\begin{claim}
If $\mu$ and $\nu$ are in $\Abb_{\Ucal}$, then $h_r(\mu \op \nu) = h_r(\mu) \op h_r(\nu)$.
\end{claim}

\begin{proof}
Let $\mu$ and $\nu$ be fixed.
Notice that, for a fixed $s$, $\nu$-a.e. $t$ satisfies
that $s << t$ and thus that
$h_r(s \op t) = h_r(s) \op h_r(t)$.
The claim now follows by unfolding the definitions:
\[
h_r(\mu \op \nu)(f) = \int \int f (x)\ d h_r(\mu \op \nu)(x)
=
\int \int f(h_r(y)) \ d\mu \op \nu (y)
\]
\[
=
\int \int f(h_r(s \op t))\ d\nu (t)\ d\mu(s)
=
\int \int f(h_r(s) \op h_r(t))\ d\nu(t)\ d\mu(s)
\]
\[
=
\int \int f(u \op v)\ d h_r(\nu)(v)\ d h_r(\mu) (u)
= h_r(\mu) \op h_r(\nu) (f)
\] 
\end{proof}

\begin{claim}
If $\mu$ is in $\Abb_{\Ucal}$, then $h_r(\mu)(E_{r,n}) = 1$ for all $n$ and
if $s \ne r$ then $\mu(E_s) = 0$.
\end{claim}

\begin{proof}
Observe that if $2^p$ divides $a+b$ and $a << b$, then
$2^p$ divides both $a$ and $b$.
Consequently, if $2^p$ divides $\#(t)$, then $h_r(t)$ is in $E_{r,2^p}$.
Now suppose that $\mu$ in $\Abb_{\Ucal}$.
Since, for every $p$, $\mu$-a.e. $t$ satisfies $2^p$ divides $\#(t)$, it
follows that $h_r(\mu) (E_{r,2^p}) = 1$.
The second conclusion of the claim follows from the observation that if
$r \ne s$, then there is an $n$ such that
$E_{r,n} \cap E_{s}$ is empty.
\end{proof}

To finish the proof of the proposition, let $\mu$ be any element of $\Abb_{\Ucal}$ such that
$\mu \op \mu = \mu$.
Define $\mu_r = h_r(\mu)$.
The above claims show that $\mu_r \op \mu_r = \mu_r$ and
that $\mu_r(E_s) = 1$ if $r = s$ and $0$ otherwise.
This finishes the proof of the proposition.
\end{proof}

In September 2012, I publicly announced a claim that Thompson's group $F$ is amenable.
The proof contained an error, found by Azer Akhmedov.
The portion of the proof that was incorrect hinged on a claim that the set of idempotent
measures in $\Abb_{\Ucal}$ was equal to a directed intersection of nonempty compact sets.
In light of the above propositions, this is impossible.

Of course the proofs in this section also suggest why this is the wrong approach: $\Abb_{\Ucal}$, which
is the closed convex hull of $\beta \Tbb \cap \Abb_{\Ucal}$, is too large.
At present, it seems plausible that if $K$ is a minimal compact subsystem of $(\beta \Tbb,\op)$,
then the closed convex hull of $K$ contains a unique idempotent.
Next we will show, however, that such $K$ are necessarily quite large.

\begin{prop}
If $K$ is a minimal compact subsystem of $(\beta \Tbb,\op)$,
then $K$ is nonseparable.
\end{prop}

\begin{proof}
If $x$ is in $2^{\Nbb}$, define $x+1$ to be the result of adding 1 to $x$ with carry to the right
(i.e. $x \mapsto x+1$ is the odometer map).
Define $h:\Tbb \to 2^{\Nbb}$ by setting $h(\one)$ to be the constant $0$ sequence and
$h(s \op t) = h(t) + 1$ and extend $h$ continuously to a map from $\beta \Tbb$ to $2^{\Nbb}$.
Notice that the $h$-image of any nonempty $\op$-closed subset of $\beta \Tbb$ is dense
and hence the image of $K$ is all of $2^{\Nbb}$ (since the odometer is a minimal dynamical system).
Moreover, for each $r$ in $2^{\Nbb}$, there is a $\xi_r$ in $K$ such that
$h(\xi_r \op \xi_r) = r$.
For each $r$ in $2^{\Nbb}$ and $p$ in $\Nbb$, define
\[
E_{r,p} = \{t \in \Tbb : \forall i < p\ (h(t)(i) = r(i))\}.
\]
Notice that $\xi_r (E_{r-1,p}) = 1$ for all $p$.
Define
\[
E_r = \{s \op t : s \op t \in E_{r,\#(s)}\} = \{s \op t : t \in E_{r-1,\#(s)} \}
\]
and observe that $\xi_r \op \xi_r (E_r) = 1$.
As in the previous section, the minimality of $K$ implies that $K \op K$
is dense in $K$.
It follows that if 
\[
U_r = \{\eta \in K : \eta(E_r) > 1/2\},
\]
then
$\{U_r : r\in 2^{\Nbb}\}$ is an uncountable family of nonempty
pairwise disjoint open subsets of $K$.
Hence $K$ is not separable.
\end{proof}

\section{Concluding remarks}

At the time this article was written, it is still unknown if $F$ is amenable.
Never-the-less, I feel 
Conjectures \ref{NA_Hindman} and \ref{NA_Ellis} are based on sound heuristics from Ramsey theory.
It is rare in the Ramsey theory of countably infinite sets that there are difficult counterexamples
to Ramsey-theoretic statements
(there are exceptions, perhaps most notably \cite{distortion}; see also
\cite[\S9]{abstract_Ramsey}).
On the other hand, there are many deep and often difficult positive results in Ramsey theory at this level:
the Dual Ramsey Theorem \cite{dual_Ramsey},
Hindman's Theorem \cite{Hindman_thm},
Gowers's $\FIN_k$ Theorem \cite{lipschitz_on_classical},
the Hales-Jewett Theorem \cite{Hales-Jewett},
and the Halpern-L\"{a}uchli Theorem \cite{Halpern-Lauchli}.
See \cite{intro_Ramsey_spaces} for further reading on these theorems as well as many others. 

Also, while we do not know whether $(\Pr(S),\star)$ contains an idempotent if $(S,\star)$ is an
arbitrary binary system, we do know that there are quite different examples
of binary systems which \emph{admit idempotent measures}:
semigroups, finite binary systems, and
binary systems depending on only one variable.

The results of this paper also suggest several test questions which allow for an incremental approach to
proving Conjectures \ref{NA_Hindman} and \ref{NA_Ellis}:

\begin{question}
Is Conjecture \ref{NA_Hindman} true for sums of $d$ elements, for a fixed $d \geq 3$?
What about the case $d =3$?
\end{question} 

\begin{question}
For which classes of binary systems is Conjecture \ref{NA_Ellis} true?
\end{question}

\begin{question}
If two binary systems satisfy Conjecture \ref{NA_Ellis}, does their product?
\end{question}

\begin{question}
For which specific values of $m$ can one prove that there is an $n$ such that
if $c:\Tbb_n \to \{0,1\}$ then there is a copy of $\Tbb_m$ in $\Abb_n$ on which $c$
is constant?
What bounds (upper or lower) can be proved on $n$ for a given value of $m$?
\end{question}

\def\Dbar{\leavevmode\lower.6ex\hbox to 0pt{\hskip-.23ex \accent"16\hss}D}

\end{document}